\documentclass[a4paper,
]{article}

\usepackage[russian,english]{babel}

\usepackage{amsthm,amsmath,amssymb,enumitem,mathrsfs}

\usepackage{graphicx}

\usepackage{scalerel}

\usepackage[utf8x]{inputenc}
\usepackage{xcolor}

\usepackage[colorlinks=true,citecolor=magenta,linkcolor=magenta]{hyperref}

\def\endproofsym{\qed}

\newenvironment{proof-nn}{\trivlist\item[\hskip\labelsep{\hskip0pt
		{\underline{Proof} :}\hskip .321429\parindent}]\ignorespaces}
{\endproofsym\endtrivlist}

\newcommand{\eps}{\varepsilon}
\newcommand{\unmezzo}{\frac{1}{2}}

\newcommand{\epsi}{\varepsilon}
\renewcommand{\epsilon}{\varepsilon}

\let\TeXchi\chi
\def\chi{{\setbox0 \hbox{\mathsurround0pt
			$\TeXchi$}\hbox{\raise\dp0 \copy0 }}}

\newcommand{\calH}{{\mathcal H}}
\newcommand{\calD}{{\mathcal D}}

\newtheorem{thm}{Theorem}[]
\newtheorem*{thm-nn}{Theorem}
\newtheorem{lem}[thm]{Lemma}
\newtheorem*{lem-nn}{Lemma}

\newtheorem*{cor-nn}{Corollary}

\newcommand{\grad}{\nabla_x} 
\newcommand{\lap}{\Delta} 

\newcommand{\R}{\mathbb{R}}

\newcommand{\sprod}[2]{\langle #1 , #2 \rangle} 

\DeclareMathOperator{\ddt}{\frac{d}{dt}}

\newcommand{\derivt}{\ddt}
\newcommand{\ints}[1]{\int_{\mathbb{T}^d}#1 \,dx} 
\newcommand{\intvp}[1]{\int_{\R^d}#1\,dv'} 

\newcommand{\intsv}[1]{\int_{\mathbb{T}^d}\int_{\R^d} #1\, dv\,dx }
\newcommand{\intsvv}[1]{\int_{\mathbb{T}^d}\int_{\R^d}\int_{\R^d} #1\, dv'\, dv\,dx }
\newcommand{\intv}[1]{\int_{\R^d} #1\, dv}

\renewcommand{\phi}{\varphi}

\def\sfT{\mathsf{T}}
\def\sfL{\mathsf{L}}

\def\sfA{\mathsf{A}}

\def\sfR{\mathsf{R}}

\newcommand{\nn}{\nonumber}

\allowdisplaybreaks

\usepackage{lipsum}

\usepackage{enumitem}
\usepackage[affil-it]{authblk}

\begin{document}


\title{\vspace{-2cm}Hypocoercivity and reaction-diffusion limit for a nonlinear generation-recombination model}

\author[1]{Gianluca Favre\footnote{E-mail: {\tt gianluca.favre@univie.ac.at}}}

\author[1]{Christian Schmeiser\footnote{E-mail: {\tt christian.schmeiser@univie.ac.at}}}

\author[2]{Marlies Pirner\footnote{E-mail: {\tt marlies.pirner@mathematik.uni-wuerzburg.de}}}

\affil[1]{Faculty of Mathematics, University of Vienna,	Oskar-Morgenstern-Platz 1, 1090 Wien, Austria}

\affil[2]{Department of mathematics, W{\"u}rzburg University, Emil Fischer Str. 40, 97074 W{\"u}rzburg, Germany}

\date{}

\maketitle

\begin{abstract}
A reaction-kinetic model for a two-species gas mixture undergoing pair generation and recombination reactions is considered on a flat torus. 
For dominant scattering with a non-moving constant-temperature background the macroscopic limit to a reaction-diffusion system is carried out. 
Exponential decay to equilibrium is proven for the kinetic model by hypocoercivity estimates. This seems to be the first rigorous derivation of a 
nonlinear reaction-diffusion system from a kinetic model as well as the first hypocoercivity result for a nonlinear kinetic problem without smallness 
assumptions. The analysis profits from uniform bounds of the solution in terms of the equilibrium velocity distribution.
\end{abstract}

\section{Introduction}
This work pursues two main directions, the derivation of a reaction-diffusion system from a kinetic reaction model and the exponential convergence of the nonlinear kinetic system towards the global equilibrium through hypocoercivity techniques. The model describes the evolution of a two species gas mixture, subject to kinetic transport, scattering with a background, a second order pair generation reaction, and its inverse recombination. As a simplification the velocity equilibria generated by the
background scattering and by the pair generation are assumed to be the same, essentially meaning that the background and the reservoir for pair generation share a
given constant temperature.
The model without the background scattering has already been
introduced in \cite{NeuSch-2016}, where a close-to-equilibrium hypocoercivity result and convergence to a nonlinear diffusion equation in the fast reaction limit have
been proven.

The results of this work provide an extension of the hypocoercivity results of \cite{NeuSch-2016}, dropping the close-to-equilibrium assumption on the initial data.
A Lyapunov functional is constructed by modification of the entropy as in \cite{DolMouSchmHypoMass}. However, instead of the quadratic close-to-equilibrium 
approximation of the entropy, the original entropy of the nonlinear problem is used. By $L^\infty$ bounds on the solution, the entropy and its dissipation are
equivalent to their quadratic approximations, which makes them compatible with the quadratic modification of the entropy introduced in \cite{DolMouSchmHypoMass}.
The hypocoercivity result also applies to the problem treated in \cite{NeuSch-2016}, since the background scattering terms are not needed for hypocoercivity, which is
a consequence of the interaction of the dissipative reaction mechanism with the mixing properties of the kinetic transport.

The formal connection between reaction-kinetic and reaction-diffusion or reaction-fluid systems is well established (see, e.g., \cite{BisiDesv2006,BisiDesv2008, SpigZan1,SpigZan2,SpigZan3,Zanette}). A rigorous justification seems to be available only for the derivation of linear reaction-diffusion systems \cite{BisiDesv2006}.
On the other hand, nonlinear diffusion and nonlinear mean field models have been derived rigorously in various contexts like semiconductors 
\cite{GolsePoupaud,PoupaudCS}, chemotaxis \cite{CMPS,HwangKangStevens}, or the fast-reaction limit in \cite{NeuSch-2016}. The compactness necessary for the
control of the nonlinear terms typically relies on the application of averaging lemmas, but also compensated compactness has been used (see, e.g., \cite{DMOS}).
Here a weighted $L^2$ averaging lemma is used, which has been proved in \cite{NeuSch-2016} in a slightly stronger version of \cite{PoupaudCS}.


The main motivation for this work is to carry out first small steps in the directions of, on the one hand, rigorously establishing the connection between nonlinear kinetic 
and macroscopic models for reacting mixtures and, on the other hand, proving quantitative decay-to-equilibrium results for nonlinear kinetic models without smallness 
assumptions. 


In the following section the model is presented together with a global existence result, based on $L^\infty$ bounds, which rely on the special structure of the reaction
mechanism. In Section 3 the nonlinear hypocoercivity result is proven giving exponential convergence to equilibrium. In the final Section 4 the convergence to a
reaction-diffusion system in the macroscopic limit is proven.

\section{The kinetic model -- global existence}
\label{sec model}
The system consists of two species with phase space distributions $f_1(x,v,t)$ and $f_2(x,v,t)$, where the position $x$ lies in the flat torus $\mathbb{T}^d$ of dimension 
$d\ge 1$ and represented by $[0,1)^d$, and the velocity $v$ lies in $\R^d$. The particles perform straight runs with velocity $v$, interrupted by
thermalizing collision events with a nonmoving medium, after which the particle velocity is sampled from an equilibrium distribution. On the other hand, pairs of particles
of both species are generated spontaneously and disappear again in recombination reactions. The simplest possible model for these generation/recombination events 
is chosen: The velocities of newly generated particles are sampled from the same equilibrium distribution as for the thermalization events, and the rate of 
the recombination reaction for one species is proportional to its phase space density and to the position density of the other species. In particular, the reaction rate is
independent from the relative velocity of the two particles. Typical applications are ionization processes or the generation/recombination of free electrons and holes
in semiconductors.

After an appropriate nondimensionalization the model takes the form 
\begin{align}
	\label{hypo1}
	\partial_t f_1 + v \cdot \grad f_1 &= \sigma \big( \rho_1 \chi_1 - f_1 \big) + \chi_1 - \rho_2 f_1 \\
	\label{hypo2}
	\partial_t f_2 + v \cdot \grad f_2 &= \sigma \big( \rho_2 \chi_2 - f_2 \big) + \chi_2 - \rho_1 f_2 \,. 
\end{align}
where the dimensionless parameter $\sigma\ge 0$ measures the relative strength of thermalization compared to generation/recombination reactions, 
and the problem treated in \cite{NeuSch-2016} is recovered with $\sigma=0$.
The position density $\rho_i$ of $f_i$ is given by
\[
    \rho_i(x,t) = \intv{f_i(x,v,t)} \,,\qquad i = 1,2\,. 
\]
For the equilibrium velocity distributions $\chi_i$ we assume that they are strictly positive, rotationally symmetric probability densities with finite fourth order
moments, i.e., 
\begin{equation}\label{ass-chi1}
   \chi_i = \chi_i(|v|^2) > 0 \,,\qquad \intv{\chi_i}=1 \,, 
\end{equation}
as well as
\begin{equation}\label{ass-chi2}
    \intv{v \chi_i}=0 \,,\qquad \intv{|v|^4 \chi_i} <\infty \,,\qquad i = 1,2\,,
\end{equation} 
with centered Gaussians as standard example.

The system \eqref{hypo1},\eqref{hypo2} will be considered subject to initial conditions
\begin{equation}\label{IC}
   f_i(x,v,0) = f_{i0}(x,v) \,,\qquad x\in\mathbb{T}^d \,,\quad v\in\R^d \,,\quad i = 1,2 \,.
\end{equation}
The essential assumption facilitating the main results of this work concerns the boundedness of the initial data, in particular we assume that 
\begin{align}
  \tag{$A_1$}\label{A1}
    \rho_m \chi_1(v) \le &\, f_{10}(x,v) \le \rho_M \chi_1(v) \,,\\
  \tag{$A_2$}\label{A2}
    \rho_M^{-1} \chi_2(v) \le &\, f_{20}(x,v) \le \rho_m^{-1} \chi_2(v) \,,
\end{align}
for $x\in\mathbb{T}^d$, $v\in\R^d$, with constants $0 < \rho_m \le \rho_M < \infty$.

Since recombination/generation happens pairwise, the difference of the total particle numbers is conserved:
\begin{equation}\label{cons-law}
\intsv{\big(f_1(x,v,t) - f_2(x,v,t)\big)} = \intsv{\big(f_{10}(x,v)-f_{20}(x,v)\big)}\,,
\end{equation}
for $t\ge 0$. The model relies on the assumption that pairs are created out of a reservoir (of neutral molecules in the case of ionization), whose density is kept constant.
A possible extension of the model would include the dynamics of this third species.


A mild formulation of the initial value problem \eqref{hypo1}, \eqref{hypo2}, \eqref{IC} is given by
\begin{eqnarray}
   f_1(x,v,t) &=& Q_2(x,v,0,t) f_{10}\left(x - vt,v\right) \label{mild1}\\
   && + \chi_1(v) \int_0^t Q_2(x,v,\tau,t) \left(1 + \sigma\rho_1\left(x + v(\tau-t),\tau\right)\right)d\tau \,,\nonumber\\
   f_2(x,v,t) &=& Q_1(x,v,0,t) f_{20}\left(x - vt,v\right) \label{mild2}\\
   && + \chi_2(v) \int_0^t Q_1(x,v,\tau,t) \left(1 + \sigma\rho_2\left(x + v(\tau-t),\tau\right)\right)d\tau \,,\nonumber
\end{eqnarray}
with
$$
   Q_i(x,v,\tau,t) = \exp\left( \sigma(\tau-t) - \int_\tau^t \rho_i\left(x+v(s-t),s\right)ds\right) \,,\qquad i = 1,2.
$$

\begin{thm}[Mild solutions]
\label{thm 1}
Let the assumptions \eqref{ass-chi1}, \eqref{ass-chi2}, \eqref{A1}, \eqref{A2} hold.
Then there exists a unique global solution of \eqref{mild1}, \eqref{mild2} such that for all $x\in\mathbb{T}^d$, $v\in\R^d$, $t\ge 0$,
\begin{align}
\rho_m \chi_1(v) \le &\, f_1(x,v,t) \le \rho_M \chi_1(v) \,,\label{est1}\\
\rho_M^{-1} \chi_2(v) \le &\, f_2(x,v,t) \le \rho_m^{-1} \chi_2(v) \,.\label{est2}
\end{align}
\end{thm}

\begin{proof-nn}
Local existence and uniqueness follows from Picard iteration applied to \eqref{mild1}, \eqref{mild2}. The iteration also propagates the
estimates \eqref{est1}, \eqref{est2}, leading to global solvability. For example, the bound $f_2 \ge \rho_M^{-1}\chi_2$, thus $\rho_2\ge \rho_M^{-1}$, 
implies 
$$
  Q_2(x,v,\tau,t)\le \exp\left( (\tau-t)\left(\sigma + \frac{1}{\rho_M}\right) \right) \,,\qquad \tau\le t \,,
$$
which, together with $f_{10}\le \rho_M\chi_1$, immediately implies the upper bound for $f_1$. The other bounds are obtained analogously.
\end{proof-nn}

\section{Hypocoercivity}


As time tends to infinity, we expect convergence of the solution of \eqref{hypo1}, \eqref{hypo2}, \eqref{IC} to a steady state $F_\infty = (f_{1,\infty},f_{2,\infty})$, 
which is spatially homogeneous, has an equilibrium velocity distribution, and where the recombination and generation terms are balanced. This 
implies 
$$
   f_{1,\infty}(v) = \rho_{1,\infty}\chi_1(v) \,,\qquad f_{2,\infty}(v) = \rho_{2,\infty}\chi_2(v) \,, 
$$
with $\rho_{1,\infty}\rho_{2,\infty} = 1$. The missing piece of information is the conservation law \eqref{cons-law}, i.e., 
$$
    \rho_{1,\infty} - \rho_{2,\infty} = \overline{\rho_{10}} - \overline{\rho_{20}} := \int_{\mathbb{T}^d} \int_{\R^d} (f_{10} - f_{20}) dv\,dx \,,
$$
implying
\begin{eqnarray*}
   \rho_{1,\infty} &=& \frac{1}{2}\left( \overline{\rho_{10}} - \overline{\rho_{20}}  \right) + \sqrt{\frac{1}{4}\left( \overline{\rho_{10}} - \overline{\rho_{20}}  \right)^2 + 1} \,,\\
   \rho_{2,\infty} &=& \frac{1}{2}\left( \overline{\rho_{20}} - \overline{\rho_{10}}  \right) + \sqrt{\frac{1}{4}\left( \overline{\rho_{10}} - \overline{\rho_{20}}  \right)^2 + 1} \,.
\end{eqnarray*}
We define the relative entropy
\begin{equation}
	\label{ln entropy}
	\calH(f_1,f_2) = \sum_{i=1,2}\intsv{\bigg[f_i \bigg( \ln\frac{f_i}{f_{i,\infty}} -1 \bigg) + f_{i,\infty} \bigg]} \,,
\end{equation}
and we compute its time derivative
\begin{align}
	\nn
	-\calD(f_1,f_2)&= - \sigma \sum_{i=1,2} \intsv{(f_i-\rho_i \chi_i) \ln \frac{f_i}{\rho_i \chi_i}}\\
	\nn
	&\quad- \intsvv{(f_1 f_2' - \chi_1 \chi_2') \ln \frac{f_1 f_2'}{\chi_1 \chi_2'}} \\
	\label{ln entropy diss}
	&=: - \sigma (\calD_1(f_1) + \calD_2(f_2)) - \calD_3(f_1,f_2) \,,
\end{align}
where obviously $\calD_i\ge 0$ for $i=1,2,3$. 

Note that $\calD = 0$ implies all properties of $F_\infty$ except spatial homogeneity, which we expect to be a consequence of the action of the conservative
transport operator. This is the essence of the hypocoercive nature of the problem. Note also that $\calD_3 = 0$ is already sufficient for these observations. The
contributions $\calD_1, \calD_2$ do not provide additional information. Therefore the analysis below also works for $\sigma=0$.

Our approach is close to the abstract method of \cite{DolMouSchmHypoMass}, which has been applied in \cite{NeuSch-2016} to the linearized problem with $\sigma=0$
in a rather straightforward way to obtain a local stability result for the nonlinear problem. Here we directly attack the nonlinear problem. As a starting point we use
the logarithmic relative entropy $\calH$ instead of its local quadratic approximation. This avoids nonlinear perturbations requiring smallness assumptions.
On the other hand, by the $L^\infty$-control the logarithmic entropy and its quadratic approximation are essentially equivalent.

For comparison with \cite{DolMouSchmHypoMass}, we write the kinetic system as an abstract ODE for $F=(f_1,f_2)^\top$ by defining the thermalization, 
the nonlinear reaction, and the transport operator, respectively, by
\begin{align*}
	\sfL F := {\rho_1 \chi_1 - f_1 \choose \rho_2 \chi_2 - f_2} \,,
	\qquad \sfR(F) := {\chi_1 - \rho_2 f_1 \choose \chi_2 - \rho_1 f_2} \,,\qquad \sfT F := v \cdot \grad F \,.
\end{align*}
Now, we can write system \eqref{hypo1}, \eqref{hypo2} as
\begin{equation}
	\label{hypo gen}
	\partial_t F + \sfT F = \sigma \, \sfL F + \sfR(F) \,.
\end{equation}
We introduce a weighted $L^2$-space appropriate for the linearized problem by defining for $F=(f_1,f_2)$ and $G=(g_1,g_2)$ the scalar product
\begin{equation*}
	\sprod{F}{G} := \sum_{i=1,2} \intsv{\frac{f_i g_i}{f_{i,\infty}}}\,,
\end{equation*}
and the induced norm $\|\cdot\|$.

The method of \cite{DolMouSchmHypoMass} would use the orthogonal projection to the null space of the linearization of the right hand side of \eqref{hypo gen}
for the definition of an auxiliary operator. Here the simpler choice
\begin{equation*}
	\Pi F = {\rho_1 \chi_1 \choose \rho_2 \chi_2} \,,
\end{equation*}
of the projection to velocity equilibrium will suffice. The auxiliary operator is then defined by
\begin{equation*}
	\sfA := \big[ 1 + (\sfT\Pi)^* (\sfT \Pi) \big]^{-1} (\sfT \Pi)^* \,,
\end{equation*}
and the modified entropy by
\begin{align}
	\nn
	\Gamma (F) &:= \calH(f_1,f_2) + \delta \sprod{\sfA (F-F_\infty)}{F-F_\infty} \\
	\label{mod entropy}
	&=\calH(f_1,f_2) + \delta \sprod{\sfA F}{F-F_\infty}
\end{align}
with $\delta >0$ to be chosen a posteriori. The second equality follows from $\sfA\Pi = 0$, a consequence of $\Pi\sfT\Pi = 0$, i.e. vanishing flux \eqref{ass-chi2}
of the velocity equilibrium, a property also used in \cite{DolMouSchmHypoMass} under the name \emph{diffusive macroscopic limit}.

We shall prove in several steps that, along solutions of \eqref{hypo1}, \eqref{hypo2}, $\Gamma$ controls the distance to the equilibrium and is controlled by
its time derivative, implying the main result of this section.

\begin{thm}[Exponential decay to equilibrium] \label{conv diss}
Let the assumptions \eqref{ass-chi1}, \eqref{ass-chi2}, \eqref{A1}, \eqref{A2} hold. Then there exist $\lambda,c>0$ such that the solution $F$ of \eqref{hypo1}, 
\eqref{hypo2}, \eqref{IC} satisfies
	\begin{equation}
		\label{exp conv F}
		\|F(\cdot,\cdot,t) - F_\infty\| \le c\, e^{-\lambda t} \|F_0 - F_\infty\| \,,\qquad t\ge 0 \,.
	\end{equation}
\end{thm}

The first step is the equivalence between the modified entropy and the norm. 

\begin{lem}
	\label{lem equi}
	Let $F$ satisfy \eqref{est1}, \eqref{est2}, and let $\delta$ be small enough. Then there exist $\alpha_2\ge\alpha_1>0$, such that
	\begin{equation*}
		\alpha_1 \|F - F_\infty \|^2 \le \Gamma (F) \le \alpha_2 \|F - F_\infty\|^2 \,.
	\end{equation*} 
\end{lem}
\begin{proof}
Taylor expansion implies
$$
   \calH (F) = \sum_{i=1}^2 \intsv{\frac{(f_i - f_{i,\infty})^2}{\hat f_i}} \,,
$$
with $\hat f_i$ taking values between $f_i$ and $f_{i,\infty}$, $i=1,2$, and therefore also satisfying \eqref{est1}, \eqref{est2}. Thus,
$$
   \frac{\rho_m}{\rho_{1,\infty}} f_{1,\infty} \le \hat f_1 \le \frac{\rho_M}{\rho_{1,\infty}} f_{1,\infty} \,,\qquad
   \frac{1}{\rho_M\rho_{2,\infty}} f_{2,\infty} \le \hat f_2 \le \frac{1}{\rho_m\rho_{2,\infty}} f_{2,\infty} \,,
$$
implying
$$
    c \|F-F_\infty\|^2 \le \calH(F) \le C \|F-F_\infty\|^2 \,,
$$
with
$$
   c = \min\left\{ \frac{\rho_{1,\infty}}{\rho_M}, \rho_{2,\infty}\rho_m \right\} \,,\qquad C = \max\left\{ \frac{\rho_{1,\infty}}{\rho_m}, \rho_{2,\infty}\rho_M \right\}\,.
$$
In \cite[Lemma 1]{DolMouSchmHypoMass} the bound 
\begin{equation}\label{A-bound}
   \|AF\| \le \frac{1}{2}\|(1-\Pi)F\|
\end{equation}
has been proven, implying
	\begin{equation*}
		\left| \sprod{\sfA (F-F_\infty)}{F-F_\infty} \right| \le \unmezzo \|F-F_\infty\|^2 \,,
	\end{equation*}
which completes the proof with $\alpha_1 := c -\delta/2$, $\alpha_2 := C + \delta/2$.
\end{proof}

The time derivative of the modified entropy is given by
\begin{eqnarray}
	\derivt \Gamma(F) &=& - \calD(F) + \delta \sprod{\sfA \sfR(F)}{F-F_\infty} + \delta \sprod{\sfA (F-F_\infty)}{\sfR(F)} \nonumber\\
	&& - \delta \sprod{\sfA \sfT F}{F-F_\infty} - \delta \sprod{\sfA F}{\sfT F} + \delta \sigma \sprod{\sfA \sfL F}{F-F_\infty} + \delta \sigma \sprod{\sfA F}{\sfL F} \nonumber\\
	&=& - \calD(F) - \delta \sprod{\sfA \sfT \Pi F}{F-F_\infty} \nonumber\\
	&&  + \delta \sprod{\sfT \sfA F}{F} - \delta \sprod{\sfA \sfT (1 - \Pi)F}{F-F_\infty} \nonumber\\
	&&+ \delta \sigma \sprod{\sfA \sfL F}{F- F_\infty} + \delta \sprod{(\sfA + \sfA^*) \sfR(F)}{F-F_\infty}\,,  \label{entr-diss}
\end{eqnarray}
where for the second equality $\sprod{\sfA F}{\sfL F}=-\sprod{(1-\Pi)\sfA F}{F} = 0$ (by $\sfA = \Pi\sfA$) has been used. It shall be proven that the first line on the right 
hand side controls the distance to equilibrium and also the terms in the second and third lines.

\begin{lem}
	\label{D_3}
	Let $F$ satisfy \eqref{est1}, \eqref{est2}. Then there exist $c_1, c_2, c_3>0$ such that 
	\begin{equation*}
		\calD_3(F) \ge c_1 \|\sfR(F)\|^2 \ge c_2 \ints{(1-\rho_1\rho_2)^2} + c_3 \|(1-\Pi)F\|^2 \,.
	\end{equation*}
\end{lem}
\begin{proof}
	We start by linearizing the logarithm in $\calD_3(F)$:
	\begin{equation*}
		\ln\frac{f_1 f_2'}{\chi_1 \chi_2'} = \frac{f_1 f_2' - \chi_1 \chi_2'}{\hat f_1 \hat f_2'} \,,
	\end{equation*}
where $\hat f_1 \hat f_2'$ takes values between $\chi_1 \chi_2'$ and $f_1 f_2'$, and therefore satisfies $\hat f_1 \hat f_2' \le \frac{\rho_M}{\rho_m} \chi_1 \chi_2'$,
implying
\begin{equation*}
	\calD_3(F) \ge \frac{\rho_m}{\rho_M} \intsvv{\frac{(f_1 f_2' - \chi_1 \chi_2')^2}{\chi_1 \chi_2'}}\,.
\end{equation*}
On the other hand, we estimate $\sfR(F)$ using the Cauchy-Schwarz inequality:
	\begin{align*}
		\|\sfR(F)\|^2 &= \sum_{i\neq j} \intsv{\frac{1}{\rho_{i,\infty} \chi_i} \bigg[\intvp{\frac{\chi_i \chi_j' - f_i f_j'}{\sqrt{\chi_j'}} \sqrt{\chi_j'} }\bigg]^2 } \\
		&\le \sum_{i\neq j} \intsv{\frac{1}{\rho_{i,\infty} \chi_i} \bigg[\bigg(\intvp{\frac{(\chi_i \chi_j' - f_i f_j')^2}{\chi_j'}}\bigg)^{1/2} \bigg(\intvp{\chi_j}\bigg)^{1/2}\bigg]^2 }\\
		&=\left(\rho_{1,\infty}+\rho_{2,\infty}\right)\intsvv{ \frac{(\chi_1 \chi_2' - f_1 f_2')^2}{\chi_1 \chi_2'}} \,.
	\end{align*}
This completes the proof of the first inequality with $c_1 = \rho_m[(\rho_{1,\infty}+\rho_{2,\infty})\rho_M]^{-1}$.
	
Alternatively, we add and subtract $\chi_i \rho_1 \rho_2$ to the $i$th component of $\sfR(F)$ and use orthogonality for computing the norm:
	\begin{align*}
		\|\sfR(F)\|^2 &= \sum_{i\neq j}\intsv{ \frac{\bigl(\chi_i (1-\rho_1 \rho_2) + \rho_j(\rho_i\chi_i - f_i)\bigr)^2}{f_{i,\infty}} }\\
		 &=  \sum_{i\neq j} \left(\frac{1}{\rho_{i,\infty}} \ints{(1-\rho_1 \rho_2)^2} + \intsv{ \rho_j^2\, \frac{(\rho_i \chi_i -f_i)^2}{f_{i,\infty}} } \right)\\
		&\ge (\rho_{1,\infty}+\rho_{2,\infty}) \ints{(1-\rho_1 \rho_2)^2}
		+ \min\left\{\rho_m^2, \rho_M^{-2}\right\} \|(1-\Pi)F\|^2 \,.
	\end{align*}
\end{proof}

\begin{lem}
	\label{lem ATP}
	With the previous definitions of $\sfT$, $\Pi$, $\sfA$, and $F_\infty$, there exists $c_4>0$ such that 
	\begin{align}
		\label{macro coer}
		\sprod{\sfA \sfT \Pi F}{F-F_\infty} &\ge c_4 \|\Pi F - \overline{\Pi} F \|^2 \,,
	\end{align}
	with 
$$
    \overline{\Pi}F := (\overline{\rho_1}\chi_1, \overline{\rho_2}\chi_2) \,,\qquad \overline{\rho_j} := \ints{\rho_j} \,,\quad j=1,2 \,.
$$
\end{lem}
\begin{proof}
The fact that $F_\infty$ is independent of $x$ implies 
$$
   \sprod{\sfA \sfT \Pi F}{F-F_\infty} = \sprod{\sfA \sfT \Pi \hat F}{\hat F} = \sprod{(1 + \mathcal{L})^{-1}\mathcal{L} \hat F}{\hat F} \,,
$$ 
with $\hat F = F - F_\infty$ and $\mathcal{L} := (\sfT \Pi)^*(\sfT\Pi)$. Obviously we need spectral information on the operator $\mathcal{L}$. Therefore we compute
$$
   \sprod{\mathcal{L}\hat F}{\hat F} = \|\sfT\Pi \hat F\|^2 = \ints{\left( \frac{D_1}{\rho_{1,\infty}} |\grad \hat\rho_1|^2 + \frac{D_2}{\rho_{2,\infty}} |\grad \hat\rho_2|^2 \right)} \,,
$$
with $\Pi \hat F = (\hat\rho_1\chi_1,\hat\rho_2\chi_2)$ and $D_j = \frac{1}{d}\intv{|v|^2 \chi_j}>0$, $j=1,2$. The Poincar\'e inequality on $\mathbb{T}^d$ with constant 
$c_P>0$ gives
\begin{eqnarray*}
  \sprod{\mathcal{L}\hat F}{\hat F} &\ge& c_P \ints{\left( \frac{D_1}{\rho_{1,\infty}} \left(\hat\rho_1 - \overline{\hat\rho_1}\right)^2
   + \frac{D_2}{\rho_{2,\infty}} \left(\hat\rho_2 - \overline{\hat\rho_2}\right)^2  \right)} \\  &\ge& \lambda_M \|\Pi \hat F - \overline{\Pi}\hat F\|^2 \,,
\end{eqnarray*}
where $\lambda_M = c_P\min\{D_1,D_2\}$ corresponds to the \emph{macroscopic coercivity} constant of \cite{DolMouSchmHypoMass}.
The application of $(1 + \mathcal{L})^{-1}\mathcal{L}$ by a spectral decomposition of $\mathcal{L}$ leads to
$$
    \sprod{\sfA \sfT \Pi F}{F-F_\infty} \ge  \frac{\lambda_M}{1 + \lambda_M} \|\Pi \hat F - \overline{\Pi}\hat F\|^2 
    = \frac{\lambda_M}{1 + \lambda_M} \|\Pi F - \overline{\Pi}F\|^2 \,,
$$
completing the proof.
\end{proof}


The previous two results can be combined to produce control of $\|F - F_\infty\|^2 = \|(1-\Pi)F\|^2 + \|\Pi F - F_\infty\|^2$:

\begin{lem}\label{lem:coercivity}
Let $F$ satisfy \eqref{est1}, \eqref{est2}. Then there exists $c_5>0$ such that
$$
    \calD_3(F) + \delta \sprod{\sfA\sfT\Pi F}{F-F_\infty} \ge c_5 \left( \|(1-\Pi)F\|^2 + \delta \|\Pi F-F_\infty\|^2\right) \,.
$$
\end{lem}

\begin{proof}
By the previous two results we have
\begin{align*}
  &\calD_3(F) + \delta \sprod{\sfA\sfT\Pi F}{F-F_\infty} \\
  &\ge c \left( \|(1-\Pi)F\|^2 + \ints{(1-\rho_1\rho_2)^2} + \delta \|\Pi F- \overline{\Pi} F\|^2\right) \,.
\end{align*}
By the conservation of the mass difference, we have $\overline{\rho_1} - \rho_{1,\infty} = \overline{\rho_2} - \rho_{2,\infty} =: r$.
Introducing $w_i = \rho_i-\overline{\rho_i}$, we have $\rho_i = \rho_{i,\infty} + r + w_i$,
$$
  \|\Pi F-\overline{\Pi} F\|^2 = \ints{\left( \frac{w_1^2}{\rho_{1,\infty}} + \frac{w_2^2}{\rho_{2,\infty}}\right)} \,,
$$
and
\begin{eqnarray*}
   \rho_1\rho_2 - 1 &=& (r + \rho_{1,\infty} + \rho_{2,\infty}) r + ( r + \rho_{2,\infty} + w_2/2)w_1 + (r + \rho_{1,\infty} + w_1/2)w_2 \\
  &=:& \gamma_0 r + \gamma_1 w_1 + \gamma_2 w_2 \,.
\end{eqnarray*}
In the right hand side, the coefficients $\gamma_1$ and $\gamma_2$ are bounded by the $L^\infty$-estimates. The coefficient of $r$, 
$\gamma_0 = \overline{\rho_1} + \rho_{2,\infty}$, is positive and bounded away from zero. This implies that the quadratic form
$$
   (\gamma_0 r + \gamma_1 w_1 + \gamma_2 w_2)^2 + \delta \left( \frac{w_1^2}{\rho_{1,\infty}} + \frac{w_2^2}{\rho_{2,\infty}}\right)
$$
is positive definite and can be estimated from below by
$$
   c\left( r^2 + \delta(w_1^2 + w_2^2) \right)
$$
with $c>0$. On the other hand, a straightforward computation gives
$$
   \|F-F_\infty\|^2 = \|(1-\Pi)F\|^2 + \left( \rho_{1,\infty} + \rho_{2,\infty}\right) r^2 
   + \ints {\left( \frac{w_1^2}{\rho_{1,\infty}} + \frac{w_2^2}{\rho_{2,\infty}}\right)} \,,
$$
completing the proof.
\end{proof}


Finally we estimate the remaining terms in the dissipation \eqref{entr-diss} of the modified entropy:

\begin{lem}
	\label{lem TAF}
	With the previous definitions of $\sfT$, $\Pi$, $\sfA$, and $F_\infty$, there exists $c_6>0$ such that 
	\begin{align}
		\label{TAF}
		|\sprod{\sfT \sfA F}{F}| &\le \|(1-\Pi)F\|^2  \,,\\
		\label{ATI-P}
		|\sprod{\sfA \sfT (1-\Pi)F}{F-F_\infty}| &\le c_6 \|(1-\Pi)F\| \sprod{\sfA \sfT \Pi F}{F- F_\infty}^{1/2} \,,\\
		\label{ALF}
		|\sprod{\sfA \sfL F}{F - F_\infty}| &\le \|(1-\Pi)F\| \sprod{\sfA \sfT \Pi F}{F- F_\infty}^{1/2} \,,\\
				\label{ARF}
		|\sprod{(\sfA + \sfA^*)\sfR(F)}{F - F_\infty}| &\le \|F - F_\infty\|  \|\sfR(F)\| \,. 
	\end{align}
\end{lem}
\begin{proof}
The estimate \eqref{TAF} has been shown in \cite[Lemma 1]{DolMouSchmHypoMass}. 

For proving \eqref{ATI-P} we follow along the lines of \cite[Section 5]{BDMMS},
and we shall use the notation from the proof of Lemma \ref{lem ATP}:
As a preparation, we note again that $\sfA \sfT\Pi F = \sfA \sfT\Pi\hat F$ since $F_\infty$ is independent of $x$.
We start by using duality and, in particular, the skew symmetry of $\sfT$:
	\begin{align} \label{AT1-P-dual}
		\sprod{\sfA\sfT(1-\Pi)F}{\hat F} &= -\sprod{(1-\Pi)F}{\sfT \sfA^* \hat F} = -\sprod{(1-\Pi)F}{\sfT^2 G} \,,
	\end{align}
with $G = \Pi (1+\mathcal{L})^{-1}\hat F = (1+\mathcal{L})^{-1}\Pi\hat F$. By its definition, $G = \Pi G$ and, thus, $G(x,v,t) = (r_1(x,t)\chi_1(v),r_2(x,t)\chi_2(v))$,
implying
$$
    \|\sfT^2 G\|^2 \le c \left(\|\nabla_x^2 r_1\|_{L^2(\mathbb{T}^d)}^2 + \|\nabla_x^2 r_2\|_{L^2(\mathbb{T}^d)}^2 \right)
    = c \left(\|\Delta_x r_1\|_{L^2(\mathbb{T}^d)}^2 + \|\Delta_x r_2\|_{L^2(\mathbb{T}^d)}^2 \right) \,,
$$
where $c$ contains bounds for the fourth order moments of $\chi_1$ and $\chi_2$. On the other hand, 
$\mathcal{L}G = (D_1\Delta_x r_1 \chi_1, D_2\Delta_x r_2 \chi_2)$ and
\begin{eqnarray}
   \sprod{\sfA \sfT \Pi \hat F}{\hat F} &=& \sprod{\mathcal{L}G}{\Pi \hat F} = \|\mathcal{L}G\|^2 + \sprod{\mathcal{L}G}{G} \ge \|\mathcal{L}G\|^2 \label{ATP-est}\\
   &=& \frac{D_1^2}{\rho_{1,\infty}} \|\Delta_x r_1\|_{L^2(\mathbb{T}^d)}^2 + \frac{D_2^2}{\rho_{2,\infty}} \|\Delta_x r_2\|_{L^2(\mathbb{T}^d)}^2
   \ge \kappa \|\sfT^2 G\|^2 \,, \nonumber
\end{eqnarray}
where we have taken the scalar product of the equation $G+\mathcal{L}G = \Pi \hat F$ with $\mathcal{L}G$. Combining this with \eqref{AT1-P-dual} gives
\eqref{ATI-P}.

By $\sfL = -(1-\Pi)$ we obtain
$$
   \sprod{\sfA \sfL F}{\hat F} = -\sprod{(1-\Pi)F}{\sfA^* \hat F} = -\sprod{(1-\Pi)F}{\sfT G} \,.
$$
As a byproduct of \eqref{ATP-est}, we have
$$
   \|\sfT G\|^2 =  \|\sfT\Pi G\|^2 = \sprod{\mathcal{L}G}{G} \le \sprod{\sfA\sfT\Pi\hat F}{\hat F} \,,
$$
completing the proof of \eqref{ALF}.

Finally, \eqref{ARF} follows from \eqref{A-bound} and from the Cauchy-Schwarz inequality.
\end{proof}

\begin{proof} (of Theorem \ref{conv diss})
We start by using Lemma \ref{lem TAF} in the entropy dissipation \eqref{entr-diss}:
\begin{eqnarray*}
	\derivt \Gamma(F) &\le& -\mathcal{D}_3(F) - \delta \sprod{\sfA\sfT F}{F-F_\infty} + \delta \|(1-\Pi)F\|^2 \\
	&&+ \delta(\sigma+c_6)\|(1-\Pi)F\| \sprod{\sfA\sfT F}{F-F_\infty}^{1/2}	+ \delta \|F-F_\infty\| \|\sfR(F)\| \,.
\end{eqnarray*}
For the further estimation we use Lemmas \ref{D_3} and \ref{lem:coercivity} and the Young inequality:
\begin{eqnarray*}
	\derivt \Gamma(F) &\le& -\frac{1}{2}\left( c_1 - \frac{\delta}{\gamma_2} \right)\|\sfR(F)\|^2 
	- \frac{\delta}{2} \left( 1 - \gamma_1(\sigma+c_6)\right)\sprod{\sfA\sfT F}{F-F_\infty} \\
	&& - \frac{1}{2}\left( c_5 - \delta \left(2 + \frac{\sigma+c_6}{\gamma_1} + \sigma\gamma_2 \right) \right) \|(1-\Pi)F\|^2 \\
	&& - \frac{\delta}{2} (c_5 - \gamma_2) \|\Pi F - F_\infty\|^2 \,.
\end{eqnarray*}
It is easily seen that, by first choosing $\gamma_1$ and $\gamma_2$ small enough, and then $\delta$ small enough, we obtain
$$
  \derivt \Gamma(F) \le -\kappa \left( \|(1-\Pi)F\|^2 + \|\Pi F - F_\infty\|^2 \right)  = -\kappa \|F-F_\infty\|^2 \,,
$$
with $\kappa > 0$. Now we employ Lemma \ref{lem equi} to obtain
$$
  \derivt \Gamma(F) \le - \lambda \Gamma(F) \,,
$$
with $\lambda = \kappa/\alpha_2$, implying $\Gamma(F) \le e^{-\lambda t}\Gamma(F_0) \le \alpha_2 e^{-\lambda t}\|F_0-F_\infty\|^2$.
Another application of Lemma \ref{lem equi} completes the proof with $c=\alpha_2/\alpha_1$.
\end{proof}

\section{Macroscopic limit}

In this section we assume that the thermalization process dominates the reactions, by setting $\sigma = \eps^{-2}$. After a rescaling of the 
position variable $x$, we obtain the system
\begin{align}
\label{eq 1}
\epsi^2 \partial_t f_1 + \epsi v \cdot \grad f_1 &= \rho_1 \chi_1 - f_1 + \epsi^2 (\chi_1 - \rho_2 f_1) \,,\\
\label{eq 2}
\epsi^2 \partial_t f_2 + \epsi v \cdot \grad f_2 &= {\rho_2} \chi_2 - f_2 + \epsi^2 (\chi_2 - \rho_1 f_2)\,,
\end{align}
which we again consider for $x\in\mathbb{T}^d$. In the following,
the macroscopic limit $\epsi\to 0$ is carried out in the initial value problem \eqref{eq 1}, \eqref{eq 2}, \eqref{IC}.
First the necessary uniform bounds are shown and then the convergence result is stated and proven. As in the previous section,
the entropy dissipation inequality
$$
  \epsi^2 \frac{d}{dt} \calH(f_1,f_2) \le -\mathcal{D}_1 - \mathcal{D}_2
$$
is used. Integration with respect to time and estimation of the dissipation terms by the squares of the norms of the weighted 
$L^2$-spaces $(L_i^2, \|\cdot\|_i)$ with 
$$
    \|f\|_i^2 := \intsv{\frac{f^2}{\chi_i}} \,,\qquad
$$
as in the previous section proves the following result.

\begin{lem}\label{lem:uniform}
Let the assumptions of Theorem \ref{thm 1} hold. Then the solution of \eqref{eq 1}, \eqref{eq 2}, \eqref{IC} satisfies 
\begin{enumerate}
\item $f_i/\chi_i$ is bounded uniformly as $\epsi\to 0$, $i = 1,2$,
\item The microscopic perturbation $R_i := (f_i -\rho_i\chi_i)/\epsi$ is bounded in $L^2((0,\infty); L_i^2)$ uniformly as $\epsi\to 0$, $i=1,2$.
\end{enumerate}
\end{lem}

Note that the lemma also implies uniform-in-$\epsi$ boundedness of the macroscopic densities $\rho_1$, $\rho_2$. Their consequential weak
convergence is not sufficient to pass to the limit in the nonlinear reaction term. Strong convergence will be obtained by employing the averaging
lemma of \cite{NeuSch-2016}. This requires a mild additional non-concentration assumption on the equilibrium distributions:
\begin{align}
\label{Ass 3}\tag{$A_3$}
&\exists\,\, C,\theta > 0: \quad \forall\,\,a \in \R, \, \omega\in \mathbb{S}^{d-1},\, \delta>0: \quad \int_{|a+v\cdot\omega|<\delta}\chi_i \,dv \le C \delta^\theta, \quad i=1,2 \,.
\end{align}

\begin{thm}[Macroscopic limit]\label{thm 2}
Let the assumptions of Theorem \ref{thm 1} and \eqref{Ass 3} hold. Then the solution of \eqref{eq 1}, \eqref{eq 2}, \eqref{IC} satisfies
$$
   \lim_{\epsi\to 0} f_i = \rho_i^0 \chi_i \qquad\mbox{in } L^2_{loc}((0,\infty);L_i^2) \,,\qquad i=1,2,
$$
where $\big(\rho_1^0(x,t),\rho_2^0(x,t)\big)$ are distributional solutions of the reaction diffusion system 
\begin{align*}
\partial_t \rho^0_1 - D_1 \lap \rho^0_1 &= 1 - \rho^0_1 \rho^0_2 \,, \\
\partial_t \rho^0_2 - D_2 \lap \rho^0_2 &= 1 - \rho^0_2 \rho^0_2 \,, 
\end{align*}
with 
$$
    D_i = \frac{1}{d} \intv{|v|^2 \chi_i} \,,\qquad i = 1,2,
$$ 
and subject to the initial conditions
$$
   \rho_i^0(t=0) = \intv{f_{i0}} \,,\qquad i=1,2\,.
$$
\end{thm}

\begin{proof-nn}
The proof follows a standard procedure for diffusive macroscopic limits. The first step is division of the equations \eqref{eq 1}, \eqref{eq 2} by $\epsi$:
\begin{align}
\epsi \partial_t f_1 + v \cdot \grad f_1 &= - R_1 + \epsi (\chi_1 - \rho_2 f_1) \,,\label{kin1/eps}\\
\epsi \partial_t f_2 + v \cdot \grad f_2 &= - R_2 + \epsi (\chi_2 - \rho_1 f_2)\,.\label{kin2/eps}
\end{align}
The uniform bounds of Lemma \ref{lem:uniform} guarantee the existence of accumulation points $\rho_i^0$ and $R_i^0$ of $\rho_i$ and, respectively,
$R_i$, $i=1,2$, as $\epsi\to 0$. We can also pass to the limit in the distributional version of the above system, to obtain
\begin{equation}\label{R0}
   R_i^0 = - \chi_i v \cdot \grad \rho_i^0 \,,\qquad i=1,2 \,.
\end{equation}
Now we divide by $\epsi$ once more and integrate with respect to $v$:
\begin{align}
  \partial_t\rho_1 + \nabla_x\cdot \intv{v R_1} &= 1 - \rho_1 \rho_2 \,,\label{macro1}\\
  \partial_t\rho_2 + \nabla_x\cdot \intv{v R_2} &= 1 - \rho_1 \rho_2 \,.\label{macro2}
\end{align}
By the Cauchy-Schwarz inequality
$$
  \left| \intv{v R_i}\right|^2 \le \intv{|v|^2 \chi_i} \intv{\frac{R_i^2}{\chi_i}}
$$
the uniform bounds of Lemma \ref{lem:uniform} suffice for passing to the limit on the left hand sides, obtaining (using \eqref{R0}) the desired
left hand sides. 

Compactness for the macroscopic densities is shown as in \cite{NeuSch-2016}: By the boundedness of the right hand sides of \eqref{kin1/eps},
\eqref{kin2/eps}, an averaging lemma \cite[Lemma 3.2]{NeuSch-2016} can be applied and gives 
$\rho_1,\rho_2 \in L^2((0,\infty);H^{\theta/(2+\theta)}(\Omega))$ uniformly in $\epsi$. The boundedness of the fluxes in the macroscopic equations
\eqref{macro1}, \eqref{macro2} implies $\rho_1,\rho_2 \in H^1((0,\infty);H^{-1}(\mathbb{T}^d))$. Combining this with the averaging lemma result by an
interpolation \cite[Lemma 3.4]{NeuSch-2016} finally gives 
$$
   \rho_1,\rho_2 \in H^{\frac{\theta}{2(1+\theta)}}((0,\infty)\times\mathbb{T}^d) \,,\qquad\mbox{uniformly in } \epsi \,.
$$
Strong convergence in $L_{loc}^2((0,\infty); L^2(\mathbb{T}^d))$ follows, permitting passage to the limit also in the right hand sides of \eqref{macro1},
\eqref{macro2}.
\end{proof-nn}

\vfill
\noindent {\bf Acknowledgments.}~~This work has been supported by the Austrian Science Fund, grants no.~W1245 and F65, and by the Humboldt foundation. G.F.~also thanks the Vienna School of Mathematics.


\label{end}

\end{document}